\documentclass{amsart}

\usepackage{amsfonts}
\usepackage[utf8]{inputenc}
\usepackage{amssymb}
\usepackage{mathabx}
\usepackage{amscd}
\usepackage{pictexwd,dcpic}
\usepackage{graphicx}
\usepackage{verbatim}
\usepackage{hyperref}
\hypersetup{
    colorlinks=true,
    linkcolor=blue,
    filecolor=magenta,
    urlcolor=cyan,
    citecolor=blue,}
\usepackage[nameinlink,capitalize]{cleveref}

\title{On the Cohen-Macaulay property of the Rees algebra of the module of differentials}
\author{Alessandra Costantini and Tan Dang}
\address{Department of Mathematics, University of California, Riverside, 900 University Ave. Riverside, CA 92521}
\email{alessanc@ucr.edu}
\address{Department of Mathematics, Purdue University, 150 N. University Street, West Lafayette, IN 47907-2067}
\email{dangt@purdue.edu}

\newtheorem{thm}{Theorem}[section]
\newtheorem{defn}[thm]{Definition}
\newtheorem{prop}[thm]{Proposition}

\newtheorem{lemma}[thm]{Lemma}
\newtheorem{cor}[thm]{Corollary}

\begin{document}
\maketitle

\begin{abstract}
   Let $R$ be an algebra essentially of finite type over a field $k$ and let $\Omega_k(R)$ be its module of K\"ahler differentials over $k$. If $R$ is a homogeneous complete intersection and $\mathrm{char}(k)=0$, we prove that $\Omega_k(R)$ is of linear type whenever its Rees algebra is Cohen-Macaulay and locally at every homogeneous prime $\mathfrak{p}$ the embedding dimension of $R_{\mathfrak{p}}$ is at most twice its dimension.
\end{abstract}

\section{Introduction}

 Let $R$ be an an algebra essentially of finite type over a field $k$. The module of K\"ahler differentials $\Omega_k(R)$ encloses relevant information about the singularities of $R$. For instance, when $R$ is reduced and local and $k$ is a perfect field, $R$ is a complete intersection if and only if $\Omega_k(R)$ has projective dimension at most one, and is regular if and only if $\Omega_k(R)$ is free \cite{{Ferrand},{Vasconcelos},{Lipman}}. Moreover, the singularities of $R$ can sometimes be understood by studying the properties of the symmetric algebra $\mathcal{S}(\Omega_k(R))$ and of the Rees algebra $\mathcal{R}(\Omega_k(R))$ \cite{{SUV1997},{SUV2012}}.

 The main goal of this article is to investigate the Cohen-Macaulay property of the Rees algebra $\mathcal{R}(\Omega_k(R))$ in the case when $\Omega_k(R)$ has projective dimension at most one. Recall that the Rees algebra $\mathcal{R}(E)$ of a module $E$ having a rank is defined as the symmetric algebra $\mathcal{S}(E)$ modulo its $R$-torsion submodule. For a module $E$ projective dimension one, the Rees algebra $\mathcal{R}(E)$ is Cohen-Macaulay whenever $E$ satisfies the so  called \emph{condition $F_1$}, a requirement on the height of the Fitting ideals of $E$ \cite{{HSV},{SUV2003}}; we refer the reader to Section~\ref{PrelimSection} for the definition of the $F_t$ condition for any integer $t \geq 0$. 

 On the other hand, modules of projective dimension one with Cohen-Macaulay Rees algebra do not necessarily satisfy condition $F_1$ (see \cite[4.7]{SUV2003}). Nevertheless, Simis, Ulrich and Vasconcelos proved the following result in the case when $E$ is the module of differentials of a complete intersection (see \cite[3.1]{SUV2012}):

 \begin{thm}[\cite{SUV2012}]\label{SUV3.1}
   Let $k$ be a field of characteristic 0 and let $R$ be a $k$-algebra essentially of finite type which is locally a complete intersection. Assume that $\Omega_k(R)$ satisfies condition $F_0$ and that $\,\mathcal{R}(\Omega_k(R))\,$ is Cohen-Macaulay. Then, $\Omega_k(R)$ satisfies condition $F_1$.
 \end{thm}

 Since $F_1$ implies $F_0$, it is natural to ask whether the conclusion of this theorem is still true if one removes the assumption that $\Omega_k(R)$ satisfies condition $F_0$ (see \cite[Question 3.3]{SUV2012}). 

 An affirmative answer to this question would imply that $\mathcal{R}(\Omega_k(R))$ is Cohen-Macaulay if and only if $\Omega_k(R)$ satisfies $F_1$. This would be a very concrete characterization, as condition $F_1$ can be easily checked through a computer algebra system. Moreover, since in this case $\,\mathcal{R}(\Omega_k(R))$ would coincide with $\,\mathcal{S}(\Omega_k(R))$ \cite{{Avramov},{HuSym},{SVsym}}, the result would imply that $\,\mathcal{R}(\Omega_k(R))$ is Cohen-Macaulay if and only if it is defined by linear equations only. Notice that this is not usually the case for modules of projective dimension one \cite[4.7]{SUV2003}.

 While not being able to remove the assumption that $R$ satisfies $F_0$, our main result, Theorem~\ref{homogeneous}, proves that in the standard graded case it suffices to assume that condition $F_0$ is satisfied only locally on the homogeneous punctured spectrum. In this case, the problem reduces to that of calculating the ideal of $t \times t$ minors $I_t(\theta)$, where $\theta$ is the transpose of the Jacobian matrix of the localization of $R$ at the homogeneous maximal ideal. This idea also appears in the proof of Theorem~\ref{SUV3.1}, where one has $t=1$. However, in our case $t \geq 2$ and the techniques used in the proof of Theorem~\ref{SUV3.1} are no longer applicable. Following a completely different path, we complete the proof by resorting to the Eagon-Northcott complex \cite{EN}. However, our methods do not apply to the local case. \\

\section{Preliminaries}\label{PrelimSection}
 In this section, we collect the necessary background on modules of differentials, as well as some known results on Rees algebras of modules of projective dimension 1 that will be needed in the proof of Theorem~\ref{homogeneous}.
 
 Throughout this paper, $k$ will be a perfect field and $R$ will denote either a finite $k$-algebra $$\,R=k[X_1, \ldots, X_n]/(f_1, \ldots, f_s),\,$$ or a $k$-algebra essentially of finite type  $$\,R=W^{-1}(k[X_1, \ldots, X_n]/(f_1, \ldots, f_s)) \,$$ for some multiplicatively closed set $W$. In both cases, the \emph{module of K\"ahler differentials} of $R$ over $k$ is the $R$-module
     $$ \Omega_k(R) \coloneq \Big( \oplus_{i=1}^{n} Rdx_i \Big) \Big/ \Big( \sum_{j =1}^{s} \frac{\partial f_j}{\partial X_i} dx_i \, \Big| \, i \in \mathcal{I} \Big) $$
    where for each $i$, $\,x_i$ is the image of $X_i$ in $R$, and $dx_i$ is the image of $x_i$ in $\Omega_k(R)$ under the universal derivation 
    $\,\displaystyle {d \colon R \longrightarrow \Omega_k(R)}$.
 In particular, one has a presentation
 $$ R^s \stackrel{\theta}{\longrightarrow} R^n \longrightarrow \Omega_k(R) \to 0 $$
 where  $\, \theta= \Big( \frac{\partial f_j}{\partial x_i} \Big)_{ij}$ is the transpose of the \emph{Jacobian matrix} of $R$. 
 
 We next list some well-known properties of the module of differentials that will be used throughout this article. We refer the reader to Kunz's book \cite{Kunz} for the proofs of these results and for an extensive treatment of modules of differentials for more general classes of rings. Recall that a finite module $E$ over a Noetherian ring has a \emph{rank} if and only if $\,E_{\mathfrak{p}} \cong R_{\mathfrak{p}}^e\,$ for every associated prime $\mathfrak{p} \in \mathrm{Ass}(R)$. In this case, $\mathrm{rank}\,E=e$. 
 
 \begin{prop} \label{basicOmega} \hypertarget{basicOmega}{}
   Suppose that $R$ is local with residue field $K= R/\mathfrak{m}$, and that the residue field extension $k \subseteq K$ is separable. Let $\mathrm{trdeg}_k K$ denote its transcendence degree. Then: 
   \begin{itemize}
       \item[(a)] The minimal number of generators of $\,\Omega_k(R)$ is 
          $$ \mu(\Omega_k(R)) = \mathrm{edim}\,R + \mathrm{trdeg}_k K, $$
        where $\mathrm{edim}\,R= \mu(\mathfrak{m})$ denotes the embedding dimension of $R$.
        \item[(b)] Assume that $R$ is reduced and $k$ is perfect. Then, $\Omega_k(R)$ has a rank if and only if $R$ is equidimensional, in which case
          $$ \mathrm{rank}(\Omega_k(R))= \mathrm{dim}\, R + \mathrm{trdeg}_k K. $$
        \item[(c)] If $\mathrm{char}(k)=0$, then $\Omega_k(R)$ has a rank if and only if $R$ is reduced and equidimensional, in which case   $\,\mathrm{rank}(\Omega_k(R))\,$ is expressed by the same formula as in part (b). 
   \end{itemize}
 \end{prop}

 Let $R$ be a Noetherian ring, and let $E$ be a finite $R$-module with $\mathrm{rank}\,E=e$. Recall that the \emph{Rees algebra} of $E$ is 
   $$ \mathcal{R}(E) \coloneq \frac{\mathcal{S}(E)}{\tau_R(\mathcal{S}(E))}\, ,$$
 where $\mathcal{S}(E)$ is the symmetric algebra and $\tau_R(\mathcal{S}(E))$ is its $R$-torsion submodule. When the symmetric algebra $\mathcal{S}(E)$ is torsion-free, the module $E$ is said to be of \emph{linear type}, since then $\mathcal{R}(E)$ is isomorphic to a quotient of a polynomial ring $R[T_1, \ldots T_{\mu(E)}]$ modulo an ideal of linear forms. When $R$ is local with maximal ideal $\mathfrak{m}$, the \emph{analytic spread} $\ell(E)$ of $E$ is defined as the Krull dimension of the \emph{special fiber ring} 
 $$\mathcal{F}(E) \coloneq \mathcal{R}(E) \otimes_R R/\mathfrak{m}\,. $$
 Moreover, the following inequalities hold (see \cite[2.3]{EHU} and \cite[2.2 and 2.3]{SUV2003}).
 
 \begin{prop} \label{analytic} \hypertarget{analytic}{}
 Let $R$ be a Noetherian local ring with dimension $d$ and let $E$ be a finite $R$-module with $\mathrm{rank}\,E=e$. Then: 
  \begin{itemize} 
      \item[(i)] If $\,d \geq 1$, then $\,e \leq \ell(E) \leq \mathrm{dim}\, R+e-1 = \mathrm{dim} (\mathcal{R}(E)) -1$.
      \item[(ii)] If the residue field of $R$ is infinite, then $\,\ell(E) \leq \mu(E)$, and equality holds if and only if $E$ has no proper reduction.
  \end{itemize}
 \end{prop}

 Recall that the $i$th \emph{Fitting ideal} $F_i(E)$ of $E$ is the ideal of minors $I_{n-i}(\varphi),\,$ where $\, \displaystyle{R^s \stackrel{\varphi}{\longrightarrow} R^n \to E \to 0}\,$ is any presentation for $E$. 
 
 \begin{defn}\label{Ft} \hypertarget{Ft}{}
    \em{Let $E$ be a finite $R$-module  with $\mathrm{rank}\,E=e$ and let $t$ be a non-negative integer. $E$ is said to satisfy \emph{condition $F_t\,$} if 
      $$ \mathrm{ht}\, F_i(E) \geq i-e+t+1 \quad \mathrm{for} \quad i \geq e\, . $$
   Equivalently, $E$ satisfies $F_t$ if and only if $\, \mu(E_{\mathfrak{p}}) \leq \mathrm{dim}\,R_{\mathfrak{p}} +e-t\, $ for every $\mathfrak{p}$ in $\mathrm{Spec}(R)$ so that $E_{\mathfrak{p}}$ is not free}.
 \end{defn}
 
 When $E$ has projective dimension one, important properties of the symmetric algebra $\mathcal{S}(E)$ and of the Rees algebra $\mathcal{R}(E)$ can be expressed via the $F_t$ conditions, as explained in the following result due to Avramov \cite[Propositions 3 and 4]{Avramov} (see also \cite[1.1]{HuSym} and \cite[3.4]{SVsym} for similar results).
 
 \begin{thm}[\cite{Avramov}] \label{projdim1} \hypertarget{projdim1}{}
   Let $R$ be a local Cohen-Macaulay ring, and let $E$ be a finite $R$-module with $\mathrm{rank}\,E=e$ and $\mathrm{projdim}\,E \leq 1$. Then: 
   \begin{itemize} 
       \item[(a)] $E$ satisfies $F_0$ if and only if $\,\mathcal{S}(E)$ is a complete intersection.
       \item[(b)] $E$ satisfies $F_1$ if and only if all the symmetric powers of $E$ are torsion-free, i.e. $E$ is of linear type.
   \end{itemize}
 \end{thm}
 
 In the situation of Theorem~\ref{projdim1}, if $E$ is $F_1$ then $\mathcal{R}(E)$ is Cohen-Macaulay \cite{{HSV},{SUV2003}}. The converse is not true in general, as shown for instance by the following result due to Simis, Ulrich and Vasconcelos (see \cite[4.7 and 4.8]{SUV2003}).
 
 \begin{thm}[\cite{SUV2003}] \label{SUV4.7} \hypertarget{SUV4.7}{}
     Let $R$ be a Gorenstein local ring with infinite residue field and let $E$ be a finite $R$-module with rank $e$ and analytic spread $\ell(E)$. Assume that $\,\mathrm{projdim}\,E=1$, that $E$ is torsion-free locally in codimension 1 and that $E$ satisfies condition $F_1$ locally in codimension $s-e+1$ for some integer $s \geq e$. Let  
     $$ 0 \to R^{n-e} \stackrel{\varphi}{\longrightarrow} R^n \to E \to 0 $$
     be a free resolution for $E$, where $n \geq s$. Then: 
       \begin{itemize} 
          \item[(a)]  $\mathcal{R}(E)$ is Cohen-Macaulay and $\,\ell(E) \leq s\,$ if and only if after elementary row operations $\, \displaystyle{I_{n-s}(\varphi)}$ is generated by the last n-s rows of $\varphi$.  
          \item[(b)] If the equivalent conditions of part (a) hold and $R$ contains an infinite field $k$, then the elementary row operations in part (a) are induced by a general element of $\mathrm{GL}_n(k)$.
      \end{itemize}
 \end{thm}

 For the module of differentials of a complete intersection ring $R$, the $F_t$ conditions can be equivalently rephrased in terms of properties of $R$. In fact, one has the following result (see \cite[2.3 and its proof]{SUV1997}).

 \begin{thm}[\cite{SUV1997}] \label{ci} \hypertarget{ci}{}
   Let $R$ be an algebra essentially of finite type over a perfect field $k$ and assume that $R_{\mathfrak{p}}$ is a reduced complete intersection for every $\mathfrak{p}$ in $\mathrm{Spec}(R)$. Then $\,\Omega_k(R)$ satisfies $F_t$ if and only if 
     $$ \mathrm{edim}\,R_{\mathfrak{p}} \leq 2\, \mathrm{dim}\, R_{\mathfrak{p}} -t $$
   for every non-regular prime $\mathfrak{p}$. \\
 \end{thm}
 
 

\section{The Rees algebra of the module of differentials} \label{MainSection}
 
 As mentioned in Section~\ref{PrelimSection}, if a module of projective dimension one satisfies $F_1$ its Rees algebra is Cohen-Macaulay. Theorem~\ref{SUV3.1} provides a partial converse in the case of the module of differentials of a complete intersection ring $R$, assuming that $\Omega_k(R)$ also satisfies $F_0$. The main goal of this paper is to understand whether the $F_0$ assumption can be removed (see \cite[Question 3.3]{SUV2012}). While not fully answering the question, our main result, Theorem~\ref{homogeneous}, proves that the $F_0$ assumption can be weakened in the case when $R$ is homogeneous.
 
 The key ingredient in the proof of Theorem~\ref{homogeneous} is Proposition~\ref{height} below. Recall that in a positively graded polynomial ring $k[X_1, \dots, X_n]$ with $\text{deg }X_i= \delta_i$, the partial derivatives of a homogeneous polynomial $f$ satisfy the \emph{Euler relation}
 $$ \sum_{1=1}^n \delta_i\, X_i\, \frac{\partial f}{\partial X_i} = (\text{deg }f)f. $$

 \begin{prop} \label{height} \hypertarget{height}{}
   Let $k$ be a field of characteristic zero, and let $R$ be a positively graded $k$-algebra $\,R=k[X_1, \dots, X_n]/I\,$ with $\mathrm{dim}\,R = d\geq 2$, where $I$ is generated by a homogeneous regular sequence $\,f_1, f_2, \dots, f_{n-d}\,$ contained in $(X_1, \dots, X_n)^2$. Consider the $n$ by $n-d$ matrix
     $$ \theta = \begin{bmatrix}
                 \frac{\partial f_1}{\partial x_1 } & \dots  & \frac{\partial f_{n-d}}{\partial x_1}\\ 
                 \vdots &  & \vdots \\ 
                 \frac{\partial f_1}{\partial x_n } & \dots  & \frac{\partial f_{n-d}}{\partial x_n}
                 \end{bmatrix} $$
  where $\,\frac{\partial f_j}{\partial x_i }$ is the image in $R$ of the partial derivative of $f_j$ with respect to $X_i$. \\
  Let $\theta'$ be the submatrix of $\theta$ consisting of the last $t$ rows where $\,t = n - 2d + 1 \geq 1$. If $\,I_t(\theta) = I_t(\theta')$, then $\,\text{ht } I_t(\theta) < d.$
 \end{prop}

 \emph{Proof}.
  Write 
  $\, \theta' = \begin{pmatrix}
               \theta_1 $\!\!\!$ & | & $\!\!\!$ \theta_2 
               \end{pmatrix} \,$ 
  where $\theta_1$ is a square $t$ by $t$ matrix and $\theta_2$ is a $t$ by $d-1$ matrix. Let $\Theta$ be the transpose of the Jacobian matrix associated with $f_1, \dots, f_{n-d}$ over the ring $k[X_1, \ldots, X_n]$. Let $\Theta'$ be the submatrix of $\Theta$ consisting of the last $t$ rows, and we write 
  $\, \Theta' = \begin{pmatrix}
             \Theta_1 $\!\!\!$ & | & $\!\!\!$ \Theta_2 
             \end{pmatrix}\,$ 
  where $\Theta_1$ is a square $t$ by $t$ matrix.

  Notice that we may reduce to the case when $\,I_1 (\Theta_2) \subset (X_1, \dots,X_{n-1}, X_n^2)$. In fact, let $\,\overline{k[X_1, \dots, X_n]} \coloneq k[X_1, \dots, X_n] / (X_1, \dots, X_{n-1})$. Since $\Theta_1$ has $t$ columns and $\Theta_2$ has $t$ rows, we can perform elementary column operations which do not change $I$, $I_t(\theta)$, and $I_t(\theta')$, until $\,I_1(\overline{\Theta_2}) \subset (\overline{X_n}^2)$. Note also that we only need to perform columns operations involving columns of the same degree, hence the entries of $\overline{\Theta}$ are still homogeneous. Thus, after replacing $\,k[X_1, \dots, X_n]\,$ with $\,\overline{k[X_1, \dots, X_n]},\,$ we may assume that $\,I_1 (\Theta_2) \subset (X_1, \dots,X_{n-1}, X_n^2)$.

  Now, let $\,\Delta_{[i_1, \dots, i_t] \times [j_1, \dots, j_t]}\,$ be the determinant of the submatrix of $\theta$ consisting of entries from rows $i_1, \dots, i_t$ and columns $j_1,\dots, j_t$. In order to simplify the notation, write $\,\Delta_{[i_1, \dots, i_t]} = \Delta_{[i_1, \dots, i_t] \times [1, \dots, t]}, \,$ $\,\partial f_i / \partial x_j = a_{ji}\,$ and $\,[1, \dots, \hat{j}, \dots, t]\,$ for the tuple of integers from $1$ to $t$ with $j$ removed. From the Laplace expansion for determinants, it follows that
  \begin{equation} \label{eq1} 
     \Delta_{[2d, \dots, n]} = \sum_{j=1}^{t} (-1)^{j+t} \, a_{nj} \, \Delta_{[2d, \dots, n-1] \times [1, \dots, \hat{j}, \dots, t]} .
  \end{equation}

 Denote $\text{deg }X_i= \delta_i$. Since all the $f_j$ are homogeneous, the Euler relations among the partial derivatives of each $f_j$ imply that in $R$ we have
  $$ \,\sum_{i=1}^{n} \delta_i x_i a_{ij}=0 \quad \mathrm{for} \: \mathrm{all} \: \,1 \leq j \leq n-d\,$$
 where for each $i$ $x_i$ is the image of $X_i$ in $R$. 
 
 Solving for $\delta_n x_n a_{nj}$ and substituting into equation (\ref{eq1}) above, we obtain that
 \begin{equation} \label{eq2}
    \begin{split}
        \delta_n x_n \,\Delta_{[2d, \dots, n]} & = \sum_{j=1}^{t} (-1)^{j+t\,} \delta_n\, x_n \,a_{nj} \,\Delta_{[2d, \dots, n-1] \times [1, \dots, \hat{j}, \dots, t]}  \\
        & = \sum_{j=1}^{t} \sum_{i=1}^{n-1}(-1)^{j+t+1\,} \delta_i \,x_i \,a_{ij} \,\Delta_{[2d, \dots, n-1] \times [1, \dots, \hat{j}, \dots, t]} \\
        & = \sum_{i=1}^{n-1\,} \delta_i \,x_i  \sum_{j=1}^t (-1)^{j+t+1\,} a_{ij} \, \Delta_{[2d, \dots, n-1] \times [1, \dots, \hat{j}, \dots, t]}\\
        & = \sum_{i=1}^{n-1}(-1)^{i\,}\delta_i \, x_i \,\Delta_{[i, 2d, \dots, n-1]\times [1, \dots, t]} ~.
    \end{split}
 \end{equation}
 Since by assumption $I_t(\theta) = I_t(\theta')$, it follows that 
  $$ \Delta_{[i, 2d, \dots, n-1] \times [1, \dots, t]} = \sum_{1 \leq \lambda_1 < \dots < \lambda_t \leq n-d} h_{i, \lambda_i, \dots, \lambda_t} \, \Delta_{[2d, \dots, n]\times [\lambda_1, \dots, \lambda_t]},$$
 for some $\, h_{i, \lambda_i, \dots, \lambda_t} \in R$. Thus we can rearrange equation (\ref{eq2}) to obtain that
 \begin{equation} \label{eq3}
   \begin{split}
     \Delta_{[2d, \dots, n]} & \Big(\delta_n \,x_n - \sum_{i=1}^{n-1}(-1)^{i\,} \delta_i \,x_i\, h_{i, 1, \dots, t} \Big)  \\ 
     & = \sum_{i=1}^{n-1}\, \sum_{[\lambda_i, \dots, \lambda_t] \ne [1, \dots, t]} (-1)^{i\,}\delta_i \,x_i \,h_{i, \lambda_1, \dots, \lambda_t\,} \Delta_{[2d, \dots, n]\times [\lambda_i, \dots, \lambda_t]}~.
   \end{split}
 \end{equation}

 In particular, equation (\ref{eq3}) provides a relation among the $t$ by $t$ minors of $\theta'$. That is, there exist elements $g_i$ and $g_i'$ in $R$ so that 
 \begin{equation} \label{eq4} 
  (\delta_nx_n - \sum_{i=1}^{n-1}x_i\, g_i, g'_2, \dots, g'_k) \in \text{ker } d_1, 
 \end{equation} 
 where $d_1$ denotes the first map of the Eagon-Northcott complex associated to $\theta'$ 
  $$ \mathbb{E}_{\bullet} \colon \,0 \rightarrow E_{n-m+1} \xrightarrow[]{} E_{n-m} \xrightarrow[]{} \dots \xrightarrow[]{d_2} E_1 \xrightarrow{d_1} E_0 \rightarrow 0 . $$

 Now, assume that $\,\text{ht } I_t(\theta')= \text{ht } I_t(\theta) = d$. Notice that this is the maximal possible height of the ideal $\, I_t(\theta'),\,$
 since $\theta'$ is a $t$ by $n-d$ matrix and $t= n-2d +1$ (see \cite[Theorem 1]{EN}). Therefore, the Eagon-Northcott complex is acyclic (see \cite[Theorem 2]{EN}). In particular, it follows that $\,\text{ker } d_1 = \text{im } d_2$. Also, by the construction of the Eagon-Northcott complex it follows that the entries in the first row of $d_2$ are in the ideal generated the entries of $\theta_2$. Hence, since $\,\text{ker } d_1 = \text{im } d_2 \,$, from equation (\ref{eq4}) it follows that
   $$\delta_nx_n - \sum_{i=1}^{n-1}x_i g_i \in I_1(\theta_2).$$
 Therefore, in $k[X_1, \dots, X_n]$ we have 
 $$ \delta_n X_n - \sum_{i=1}^{n-1}X_i G_i \in I_1(\Theta_2) + I $$ 
 where for each $i$ $G_i$ denotes the preimage of $g_i$ in $k[X_1, \dots, X_n].$ 

 But this is impossible because $\delta_n$ is a unit, $\,I_1(\Theta_2) \subset (X_1, \dots,X_{n-1}, X_n^2)\,$ and $\,I \subset (X_1, \dots, X_n)^2$.
 $\blacksquare$\\

 The following lemma appeared in a slightly less general setting in the proof of Theorem~\ref{SUV3.1} (see \cite[proof of 3.1]{SUV2012}).
 
 \begin{lemma} \label{algebraic} \hypertarget{algebraic}{}
   Let $k$ be an infinite perfect field, and let $R$ be a reduced local Cohen-Macaulay $k$-algebra essentially of finite type with $\,\mathrm{dim}\,R \geq 1$. Assume that the residue field $K$ of $R$ is infinite and that $\,\mathcal{R}(\Omega_k(R))\,$ is Cohen-Macaulay. Then, there exists a finitely generated field extension $L$ of $K$ so that $L$ is algebraic over $k$ and the Rees algebra $\,\mathcal{R}(\Omega_L(R))\,$ is Cohen-Macaulay.
 \end{lemma}
 
 \emph{Proof}.
    Let $t= \mathrm{trdeg}_k(K)$ and suppose that $t \geq 1$. Write $\,R=(k[x_1, \ldots, x_m])_{\mathfrak{p}}\,$ for some prime ideal $\mathfrak{p}$ of $k[x_1, \ldots, x_m]$, and let $\, y_1, \ldots, y_t\,$ be a transcendence basis of $K$ over $k$. Since $k$ is infinite, we may assume that the $y_i$'s are general $k$-linear combinations of $x_1, \ldots, x_m$. Let $L \coloneq k(y_1, \ldots, y_t)$. Then, $L$ is a field with $k \subseteq L \subseteq R$, hence the Zariski sequence (see for instance \cite[Proposition 16.2]{Eis})
      $$ R \otimes \Omega_k(L) \longrightarrow \Omega_k(R) \longrightarrow \Omega_L(R) \to 0 $$
    implies that 
       $$ \Omega_L(R) \cong \frac{\Omega_k(R)}{R dy_1+ \ldots +R dy_t} $$
    where $d \colon R \longrightarrow \Omega_k(R)$ is the universal derivation. 
    
    On the other hand, by Proposition~\ref{basicOmega}(b) one has that $$\,\mathrm{rank}(\Omega_k(R))= t+ \mathrm{dim}\,R \geq 2.$$ 
    Hence, since $\mathcal{R}(\Omega_k(R))$ is Cohen-Macaulay and $y_1, \ldots, y_t$ are general elements, by \cite[2.2(f)]{CPU} it follows that 
       $$ \mathcal{R}(\Omega_L(R)) \cong \frac{\mathcal{R}(\Omega_k(R))}{(dy_1, \ldots , dy_t)}, $$
    so $\,\mathcal{R}(\Omega_L(R))\,$ is Cohen-Macaulay. $\blacksquare$\\

 We are finally ready to prove our main result.
 
 \begin{thm}\label{homogeneous} \hypertarget{homogeneous}{}
   Let $k$ be a field of characteristic zero. Let $R = k[X_1, \ldots X_n]/I$ be a standard graded $k$-algebra where $I$ is generated by a homogeneous regular sequence. Assume that $d \coloneq \mathrm{dim}\,R \geq 1$ and that the following conditions hold.
   \begin{itemize}
     \item[(i)] $\mathrm{edim}\, R_{\mathfrak{p}} \leq 2 \, \mathrm{dim}\, R_{\mathfrak{p}}$ for every prime $\mathfrak{p}$ in $\mathrm{Spec}(R)\setminus V(R_{+})$;
     \item[(ii)] $\mathcal{R}(\Omega_k(R))$ is Cohen-Macaulay.
  \end{itemize}
     Then $\,\mathrm{edim}\, R_{\mathfrak{p}} \leq 2 \, \mathrm{dim}\, R_{\mathfrak{p}}-1\,$ for every non-minimal homogeneous prime $\mathfrak{p}$ in $\mathrm{Spec}(R)$.
 \end{thm}

 \emph{Proof}. Let $\mathfrak{m}=(X_1, \ldots, X_n)$ be the unique homogeneous maximal ideal, then we may assume that $I \subseteq \mathfrak{m}^2$. Let $\,f_1, \ldots, f_{n-d}\,$ be a homogeneous regular sequence generating $I$. Notice that assumption (i) implies that $R$ is reduced. Hence, by Proposition~\ref{basicOmega}(c) $\,\Omega_k(R)$ has a rank, and $\mathrm{rank}(\Omega_k(R)) \geq d$. Moreover, by Theorem~\ref{SUV3.1} (see \cite[3.1]{SUV2012}) it follows that $\,\mathrm{edim}\, R_{\mathfrak{p}} \leq 2 \, \mathrm{dim}\, R_{\mathfrak{p}}-1\,$ for every non-minimal homogeneous prime $\mathfrak{p}$ with $\mathrm{dim}R_{\mathfrak{p}} \leq d-1$. Thus, it only remains to show that $\, n=\mathrm{edim}\,R_{\mathfrak{m}} \leq 2\, \mathrm{dim}\,R_{\mathfrak{m}}-1=2d-1$.
 
 After localizing at $\mathfrak{m}$, we may assume that $R$ is a reduced local complete intersection and that the Rees algebra $\mathcal{R}(\Omega_k(R))$ is Cohen-Macaulay (recall that formation of Rees algebras commutes with flat base change, see \cite[1.3]{EHU}). If $d=1$, since $\,\mathcal{R}(\Omega_k(R))\,$ is Cohen-Macaulay, by \cite[4.3]{SUV2003} it follows that $\Omega_k(R)$ modulo its $R$-torsion submodule is free. By the Jacobian Criterion this means that $R$ is regular (because $\mathrm{char}(k)=0$, see \cite[Theorem]{Lipman}). Hence, $n=d=1=2d-1$ and the statement is true in this case.
 
 Now, let $d\geq 2$ and suppose by contradiction that $\, n \geq 2d$. Then, $\mathrm{rank}(\Omega_k(R)) \geq 2$. Since the Rees algebra is Cohen-Macaulay, by Lemma~\ref{algebraic} there exists a finitely generated field extension $L$ of $K$ so that $L$ is algebraic over $k$ and the Rees algebra $\mathcal{R}(\Omega_L(R))$ is Cohen-Macaulay. Hence, after replacing $L$ with $k$ we may assume that the residue field $K$ of $R$ is algebraic over $k$. Thus, $\, \mathrm{rank}(\Omega_k(R))= d$ and $\, \mu(\Omega_k(R)) =n$ by Proposition~\ref{basicOmega}. Hence, $\,\Omega_k(R)$ is a module of projective dimension one, with a minimal free resolution
    $$ 0 \to R^{n-d} \stackrel{\theta}{\longrightarrow} R^n \longrightarrow \Omega_k(R) \to 0 $$
 where $\,\displaystyle{\theta= \Big( \frac{\partial f_j}{ \partial x_i} \Big)}\,$ and for each $i$, $\frac{\partial f_j}{ \partial x_i}$ denotes the image of the partial derivative of $f_j$ with respect to $X_i$ in $R$.
 
 Notice that $\,\Omega_k(R)\,$ has analytic spread $\, \ell(\Omega_k(R)) \leq d+ \mathrm{rank}(\Omega_k(R)) -1 =2d-1$. Since $\Omega_k(R)$ satisfies $F_1$ locally in codimension $d-1$ and $\mathcal{R}(\Omega_k(R))$ is Cohen-Macaulay, we can then apply  Theorem~\ref{SUV4.7} with $\,s=2d-1 \leq n-1\,$ to deduce that, after elementary row operations over $k$, $I_t(\theta) = I_t(\theta')$ where $\theta'$ is the submatrix of $\theta$ consisting of the last $t$ rows, with $\, t = n  - 2d + 1$. By Proposition~\ref{height} it follows that $\,\mathrm{ht }\, I_t(\theta) < d$. However, by assumption (i), locally at any prime $\mathfrak{q}$ containing $I_t(\theta)$, one has $\, \mathrm{ht}(I_t(\theta)_{\mathfrak{q}}) \geq n-d-t+1 =d$, which is a contradiction. $\blacksquare$\\

 \begin{cor}\label{corhomogeneous}\hyperlink{corhomogeneous}{}
   Let $k$ be a field of characteristic zero. Let $R = k[X_1, \ldots X_n]/I$ be a standard graded $k$-algebra where $I$ is generated by a homogeneous regular sequence. Assume that $\,\mathrm{edim}\, R_{\mathfrak{p}} \leq 2 \, \mathrm{dim}\, R_{\mathfrak{p}}\,$ for every prime $\mathfrak{p}$ in $\mathrm{Spec}(R)\setminus V(R_{+})$.
   
   Then $\,\mathcal{R}(\Omega_k(R))$ is Cohen-Macaulay if and only if $\mathrm{edim}\, R_{\mathfrak{p}} \leq 2 \, \mathrm{dim}\, R_{\mathfrak{p}}-1$ for every non-minimal homogeneous prime $\mathfrak{p}$ in $\mathrm{Spec}(R)$.
 \end{cor}

 \emph{Proof}. The forward direction follows from Theorem~\ref{homogeneous}. Now, notice that since $R$ is standard graded condition $F_1$ is satisfied if and only if it is satisfied locally at homogeneous primes. Hence, the backward direction follows from Theorem~\ref{projdim1} and Theorem~\ref{ci}. $\blacksquare$\\
 
 Theorem~\ref{homogeneous} shows that the Cohen-Macaulay property of the Rees algebra puts strong geometric constraints on the module of differentials. In fact, one has the following result.

 \begin{cor}
   Let $k$ be a field of characteristic zero. Let $V \subset \mathbb{P}^n$ be a smooth non-degenerate complete intersection subvariety of dimension $d$. Let $R$ be the homogeneous coordinate ring of $V$. 

   Then $\,\mathcal{R}(\Omega_{R/k})$ is Cohen-Macaulay if and only if $n \leq 2d$.
 \end{cor}

 \emph{Proof}. The statement follows immediately from Corollary~\ref{corhomogeneous} after we notice that $\text{edim }R = n+1$ and we only need to check the inequality at the unique homogeneous maximal ideal. $\blacksquare$\\

  Notice that the proof of Theorem~\ref{homogeneous} only requires the assumption that $R$ is standard graded in order to perform the calculation of the height of the ideal $\, I_t(\theta)$ via Proposition~\ref{height}, in the case where the dimension of $R$ is at least 2. This calculation was performed using an entirely different technique in the case of Theorem~\ref{SUV3.1} where $t=1$ (see proof of \cite[3.1]{SUV2012}). Hence, it is natural to ask whether alternative methods to calculate $\, \mathrm{ht}\,I_t(\theta)$ could make it possible to extend Theorem~\ref{homogeneous} to more general cases.

\section*{Acknowledgements}
   This problem was suggested to the authors by their Ph.D. advisor Prof. Bernd Ulrich. The authors are deeply grateful to Professor Ulrich for many interesting conversations and insightful comments on this topic.\\

\end{document}